\documentclass[11pt,reqno]{amsart}

\usepackage{preamble}

\hypersetup{colorlinks=true,
    linkcolor=blue,
    citecolor=blue,
    pdftitle={Distance sets bounds for polyhedral norms via effective dimension},
    pdfauthor={Iqra Altaf, Ryan Bushling, and Bobby Wilson}}

\begin{document}

\begin{abstract}
    We prove that, for every norm on $\mathbb{R}^d$ and every $E \subseteq \mathbb{R}^d$, the Hausdorff dimension of the distance set of $E$ with respect to that norm is at least $\hdim E - (d-1)$. An explicit construction follows, demonstrating that this bound is sharp for every polyhedral norm on $\mathbb{R}^d$. The techniques of algorithmic complexity theory underlie both the computations and the construction.
\end{abstract}

\title[Distance sets bounds via effective dimension]{Distance sets bounds for polyhedral norms \\ via effective dimension}
\author[Altaf]{Iqra Altaf}
\address{Department of Mathematics \\ University of Chicago \\ 5734 S University Avenue, Room 108 \\ Chicago, IL, 60637-1505}
\email{iqra@uchicago.edu}
\author[Bushling]{Ryan Bushling}
\address{Department of Mathematics \\ University of Washington, Box 354350 \\ Seattle, WA 98195-4350}
\email{reb28@uw.edu}
\author[Wilson]{Bobby Wilson}
\address{Department of Mathematics \\ University of Washington, Box 354350 \\ Seattle, WA 98195-4350}
\email{blwilson@uw.edu}
\thanks{B.~Wilson was supported by NSF grant, DMS 1856124, and NSF CAREER Fellowship, DMS 2142064.}


\maketitle

\section{Introduction} \label{s:intro}

\noindent Given a set $E \subseteq \R^d$, let
\begin{equation*}
    \Delta(E) := \{ \| x-y \| \in [0,\infty) :~ x,y \in E \}
\end{equation*}
be the \textit{distance set} of $E$, where $\| \cdot \|$ is the Euclidean norm. The \textit{Falconer distance problem} is to find a lower bound on $\hdim \Delta(E)$, the Hausdorff dimension of the distance set, in terms of the Hausdorff dimension of $E$ itself. A frequent simplification is to fix a point $x \in \R^d$ and instead investigate the \textit{pinned distance set}
\begin{equation*}
    \Delta_x(E) := \{ \| x-y \| \in [0,\infty) :~ y \in E \}.
\end{equation*}
Since $\Delta_x(E) \subseteq \Delta(E)$ whenever $x \in E$, any lower bound on the size of $\Delta_x(E)$ immediately implies the same bound on the size of $\Delta(E)$. An interesting variant of the Falconer distance problem is to seek dimension bounds on distance sets relative to a different norm $\| \cdot \|_*$. In this case, the corresponding distance and pinned distance sets are denoted $\Delta^*(E)$ and $\Delta_x^*(E)$, respectively.

Falconer \cite{falconer1985hausdorff} proves a simple estimate in the Euclidean case:
\begin{thm}[Falconer \cite{falconer1985hausdorff}]
    If $E \subseteq \R^d$ is any set, then 
    \begin{equation*}
        \hdim \Delta(E) \geq \hdim E - (d-1).
    \end{equation*}
\end{thm}
The elementary proof presented by Falconer generalizes to other norms with little modification. However, the techniques of algorithmic complexity provided in \cite{lutz2018algorithmic} are highly applicable here, and we use this theory to give an alternate proof for pinned distance sets of the general case.
\begin{thm} \label{thm:bound}
    Let $\| \cdot \|_*$ be any norm on $\R^d$ and let $E \subseteq \R^d$. Then
    \begin{equation*}
        \hdim \Delta_x^*(E) \geq \hdim E - (d-1)
    \end{equation*}
    for all $x \in \R^d$.
\end{thm}

\begin{rem}
    The proof of Theorem \ref{thm:bound} does not use the symmetry of $\| \cdot \|_*$, so the theorem statement can be generalized to so-called \textit{asymmetric norms}. The same is true of Theorem \ref{thm:example} below.
\end{rem}

Polyhedral norms provide a class of finite-dimensional normed spaces which can be considered the ``worst-case scenario" for the distance set problem, as the curvature of the norm ball plays a major role in precluding the possibility of a large number of points with equal pairwise distances. The relevance of nonvanishing curvature manifests itself in \cite{guth2020falconer}: the authors provide the best known bound in the planar Euclidean case, and it is demonstrated that the essential characteristic of the Euclidean ball required for their argument is that it is a smooth norm ball with strictly positive curvature. A construction of Falconer \cite{falconer2004dimensions} illustrates the marked differences in the context of polyhedral norms.
\begin{prop}[Falconer \cite{falconer2004dimensions}]
    Let $\| \cdot \|_P$ be a polyhedral norm on $\R^d$. Then there exists a compact set $F \subset \R^d$ with $\hdim F = d$ such that $\mathcal{L}(\Delta^P(F))$ $= 0$.
\end{prop}
The relevance of curvature to the distance set conjecture is further explored in \cite{bishop2021falconer}, in which Falconer's construction is generalized to the case in which at most countably many points of the surface of the unit ball do not lie on any open line segment contained in the surface of the unit ball (all but finitely many points on the surface of a polyhedral ball satisfy this condition) or for certain convex sets for which the measure representing the curvature field is singular with respect to surface measure. For more state-of-the-art estimates and information on the history of the Euclidean distance set conjecture, see \cite{du2021improved} and \cite{du2023new}.

 Our second result shows that polyhedral norms witness the sharpness of the dimension bound in Theorem \ref{thm:bound}. This, too, will be proven using complexity theoretic tools.

\begin{thm} \label{thm:example}
    Let $\| \cdot \|_P$ be a polyhedral norm on $\R^d$. For any $s \in [d-1,d]$, there exists a compact set $E \subset \R^d$ with $\hdim E = s$ such that
    \begin{equation} \label{eq:sharp}
        \hdim \Delta^P(E) = s - (d-1).
    \end{equation}
\end{thm}

\begin{rem}
    The proof of Theorem \ref{thm:bound} also works with \textit{Packing dimension} in place of Hausdorff dimension. However, the resulting statement for Packing dimension is far from sharp.
\end{rem}

Under the hypotheses of Theorem \ref{thm:example}, it follows from Theorem \ref{thm:bound} that the pinned distance sets of $E$ have minimal dimension $s - (d-1)$ at every pin in $E$: that is,
\begin{equation*}
    \hdim \Delta_x^P(E) = s - (d-1) \quad \text{for all} \quad x \in E.
\end{equation*}

Following a brief synopsis of the background material, we prove Theorem \ref{thm:bound} in \S \ref{s:bound}. In \S \ref{s:examples} we construct the examples described in Theorem \ref{thm:example}.

\section{Preliminaries} \label{s:prelim}

This section summarizes the requisite terminology and results from algorithmic complexity theory. See \cite{lutz2018algorithmic} for a more thorough exposition of this material. We denote by
\begin{equation*}
    \{0,1\}^* := \bigcup_{n \in \N} \{0,1\}^n
\end{equation*}
the set of all (finite) binary strings, including the empty string $\lambda \in \{0,1\}^0$.

\begin{defn}[Kolmogorov Complexity of a String]
    Let $\sigma, \tau \in \{0,1\}^*$. The \textit{conditional Kolmogorov complexity of $\sigma$ given $\tau$} is the length of the shortest program $\pi$ that will output $\sigma$ given $\tau$. More precisely,
        \begin{equation*}
            K(\sigma|\tau) := \min_{\pi \in \{0,1\}^*} \{ \ell(\pi) :~ U(\pi,\tau) = \sigma\},
        \end{equation*}
    where $U$ is a fixed universal prefix-free Turing machine and $\ell(\pi)$ is the length of $\pi$. The \textit{Kolmogorov complexity of $\sigma$} is simply the conditional Kolmogorov complexity of $\sigma$ given the empty string:
    \begin{equation*}
        K(\sigma) := K(\sigma|\lambda),
    \end{equation*}
    where $\lambda$ is the empty string.
\end{defn}

Throughout, we work with a fixed \textit{encoding} $\{0,1\}^* \to \bigcup_{d \in \N} \Q^d$, under which the definitions above extend from strings to vectors over the rationals.

\begin{defn}[Kolmogorov Complexity of a Point]
    Let $x \in \R^d$ and $r \in \Z_+$. The \textit{Kolmogorov complexity of $x$ at precision $r$} is the length of the shortest program that outputs a point in $\Q^d$ that approximates $x$ to $r$ bits of precision:
    \begin{equation*}
        K_r(x) := \min \, \{ K(p) :~ p \in B(x,2^{-r}) \cap \Q^d \}.
    \end{equation*}
    If also $y \in \R^{d'}$ and $s \in \Z_+$, then the \textit{conditional Kolmogorov complexity of $x$ at precision $r$ given $y$ at precision $s$} is defined by
    \begin{equation*}
        K_{r,s}(x|y) := \max \, \big\{ \min \{ K(p|q): p \in B(x,2^{-r}) \cap \Q^d \}: q \in B(y,2^{-s}) \cap \Q^{d'} \big\}.
    \end{equation*}
\end{defn}

Note that, for $x \in \R^d$, $K_r(x)$ is always less than equal to $dr + O(\log r)$. To simplify notation, we also denote $K_{r,r}(x|y)$ by $K_r(x|y)$, $K_{r,s}(x|x)$ by $K_{r,s}(x)$, and, when $r \in (0,\infty)$, $K_{\lceil r \rceil}(x)$ by $K_r(x)$.

\begin{lem}[Symmetry of Information \cite{lutz2020bounding}] \label{lem:symmetry}
    For all $m,n \in \N$, $x \in \R^m$, $y \in \R^n$, and $r,s \in \N$ with $r \geq s$: \vs{-0.15}
    \begin{enumerate}[label={\normalfont\textbf{\alph*.}},itemsep=3pt,topsep=0pt]
        \item $\big| K_r(x|y) + K_r(y) - K_r(x,y) \big| \leq O(\log r) + O(\log \log \|y\|)$.
        \item $\big| K_{r,s}(x|x) + K_s(x) - K_r(x) \big| \leq O(\log r) + O(\log \log \|x\|)$.
    \end{enumerate}
\end{lem}

In the remainder of the article, we will use the $O(\log r)$ error term to encapsulate all errors accumulated from computing terms that do not change as the precision, $r$, increases to infinity. This would include the $\log \log \|y\|$ and $\log \log \|x\|$ terms from Lemma \ref{lem:symmetry} as well as $\log m$ and $\log n$ terms already suppressed in the statement of Lemma \ref{lem:symmetry}. Using the notion of Kolmogorov complexity, one defines the concept of effective dimension as the asymptotic Kolmogorov complexity of a given point.
\begin{defn}
    The \textit{effective Hausdorff dimension} of a point $x \in \R^d$ is given by
    \begin{equation*}
        \dim(x) := \liminf_{r \to \infty} \frac{K_r(x)}{r}.
    \end{equation*}
The \textit{effective Packing dimension} of a point $x \in \R^d$ is given by
    \begin{equation*}
        \Dim(x) := \limsup_{r \to \infty} \frac{K_r(x)}{r}.
    \end{equation*}
\end{defn}

In computability theory, a set $A \subseteq \N$ is called an \textit{oracle}. Heuristically, a Turing machine $T$ can be said to \textit{have access to $A$} if, in addition to its usual operations, it can inquire whether the number currently printed on the work tape belongs to $A$. This operation takes one step, and the answer to the question is recorded as a state. The machine derived from $T$ by allowing it access to $A$ is denoted $T^A$.

The concepts of complexity and dimension defined above can be \textit{relativized to an oracle $A$} by replacing the fixed universal Turing machine $U$ with $U^A$. These relativized concepts are denoted $K^A$, $\dim^A$, etc. For us, the utility of oracles is that they allow us to compute the Hausdorff dimension of a set from the effective dimensions the points it contains.
\begin{thm}[Point-to-Set Principle \cite{lutz2018algorithmic}] \label{thm:pts}
    For every set $E \subseteq \R^d$,
    \begin{equation*}
        \hdim E = \min_{A \subseteq \N} \sup_{x \in E} \dim^A(x) \quad \text{and}\quad \pdim E = \min_{A \subseteq \N} \sup_{x \in E} \Dim^A(x).
    \end{equation*}
\end{thm}
One can infer from the proof of the point-to-set principle that, for any oracle $A$ and $s \in (0,d]$, the set of points $x$ that satisfy $\dim^A(x) < s$ has $s$-dimensional Hausdorff measure zero. In particular, $\dim^A(x) = d$ for all $x$ outside of a Lebesgue-null set.

As in \cite{stull2022pinned}, we describe the complexity of a point $x \in \R^d$ relative to a point $y \in \R^{d'}$ as the complexity of $x$ relative to an oracle set $A_y$ that encodes the binary expansion of $y$ in a standard way. We define 
\begin{equation*}
    K_r^y(x) := K_r^{A_y}(x).
\end{equation*}
Throughout the paper, we will often invoke the notion of (Martin--L\"{o}f) ``randomness." A point $x \in \R^d$ is \textit{random} with respect to an oracle $A$ if
\begin{equation*}
    \big|K_r^A(x) - rd\big| \leq O(\log r). 
\end{equation*}
This also implies that each of the coordinates of $x$ is random with respect to the others and has effective Hausdorff dimension $1$ with respect to the oracle $A$. More quantitatively, for $t < r$, the symmetry of information (Lemma \ref{lem:symmetry}) implies that
\begin{align*}
    rd + O(\log r) = K_r^A(x) &= K_{r,t}^A(x) + K_t(x) + O(\log r) \\
    &\leq K_{r,t}^A(x) + td + O(\log r),
\end{align*} 
whence
\begin{equation*}
    |K_{r,t}^A(x) - (r-t)d| \leq O(\log r).
\end{equation*}
This means that the first $t$ digits of the binary expansion of $x$ ``do not help much" in the calculating the remaining $r-t$ digits, and that we would approximately need $(r-t)d$ bits to calculate them. In the proof of Theorem \ref{thm:example}, we will have occasion to choose collections of strings of digits that are random with respect to each other, but we can see from this discussion that such a choice is always possible.

\section{Distance set bounds for arbitrary norms} \label{s:bound}

We begin with a ``trivial" lemma for the benefit of the reader new to the theory of Kolmogorov complexity.

\begin{lem}
    Let $A \subseteq \N$. For all $x,y \in \R^d$ and $r \in \N$,
    \begin{equation} \label{eq:addition}
        K_r^{A,x}(y) \leq K_r^{A,x}(x+y) + O(\log r).
    \end{equation}
\end{lem}

\textit{Proof.} Let $T^{A,x}$ be an oracle Turing machine that operates as follows: given a rational point $p \in \Q^d$, a precision level $r \in \N$, and a ``corrector" string $\sigma = \sigma_1 \dots \sigma_{d'} \in \{0,1\}^{d'}$ of length $d' := \big\lceil d^{1/2} \big\rceil$, it first computes the ``small" dyadic rational number
\begin{equation*}
    q_{\sigma,r} := \frac{\sigma_1}{2^{-(r+1)}} + \cdots + \frac{\sigma_{d'}}{2^{-(r+d')}}.
\end{equation*}
The machine then returns the first $r+d'+1$ bits of $p + q_{r,\sigma} - x$ as its output:
\begin{equation*}
    T^{A,x}(p,r,\sigma) = (p + q_{r,\sigma} - x) \upharpoonright (r+d'+1).
\end{equation*}
Because the oracle encodes $x$, this is a computable operation relative to $(A,x)$; hence, $T^{A,x}$ exists.

Now let $p \in B(x+y,2^{-r}) \cap \Q^d$ be a point with $K^{A,x}(p) = K_r^{A,x}(x+y)$. We use $p$ and $T^{A,x}$ to compute a precision-$r$ approximation of $y$. To do this, take $\sigma \in \{0,1\}^{d'}$ such that $p + q_{r,\sigma} \in B(x+y,2^{-(r+1)})$. (This explains the choice of $d'$: it is just large enough to encode a rational number that ``corrects" $p$ by a factor of up to $\tfrac{1}{2}$, upgrading the precision-$r$ approximation to a precision-$(r+1)$ approximation.) Since $\| (z \upharpoonright (t+d')) - z \| < 2^{-t}$ for all $z \in \R^d$ and $t \in \N$,
\begin{align*}
    &\big\| \big( (p + q_{r,\sigma}) - x) \upharpoonright (r+d'+1) \big) - y \big\| \\
    &\qquad \leq \big\| (p + q_{r,\sigma} - x) - y \big\| + \big\| \big( (p + q_{r,\sigma} - x) \upharpoonright (r+d'+1) \big) - (p + q_{r,\sigma} - x) \big\| \\
    &\qquad < \| (p + q_{r,\sigma}) - (x+y) \| + 2^{-r-1} < 2^{-(r+1)} + 2^{-(r+1)} = 2^{-r},
\end{align*}
so $T^{A,x}(p,r,\sigma) \in B(y,2^{-r})$.

If $\mu$ is a string that encodes $T^{A,x}$ in the language of $U^{A,x}$, so that $U^{A,x}(\mu,q,r,\sigma) = T^{A,x}(q,r,\sigma)$ for all $q$, $r$, and $\sigma$, then 
\begin{align*}
    K_r^{A,x}(y) &\leq K^{A,x}(\mu,p,r,\sigma) \leq |\mu| + K^{A,x}(p) + K^{A,x}(r) + K^{A,x}(\sigma) + O(1) \\
    &= K^{A,x}(p) + O(\log r) = K_r^{A,x}(x+y) + O(\log r).
\end{align*}
(The value $|\mu|$ is a constant independent of $r$, and the additional $O(1)$ represents the complexity required to concatenate $\mu$, $p$, $r$, and $\sigma$ into a single string in the language of $U^{A,x}$.) \textqed

This computations of this sort are routine, arithmetic properties such as that in \eqref{eq:addition} will be used freely in what follows.

\noindent \define{Proof of Theorem \ref{thm:bound}.} The unit sphere in the norm $\| \cdot \|_*$ has Packing dimension $d-1$. Therefore, by the point-to-set principle, there exists an oracle $A \subseteq \N$ such that
\begin{equation*}
    \sup_{z \in \partial B_*(0,1)} \Dim^A(z) = d-1.
\end{equation*}
For any  $y \in E$, if we write $y = x + \|x-y\|_* z$ for some $z \in \partial B_*(0,1)$, then the following holds: for any oracle $B \subseteq \N$, 
\begin{align*}
    K_r^{A,B,x}(y) & \leq K_r^{A,B,x}(y-x) + O(\log r) \\
    & = K_r^{A,B,x}(\|x-y\|_* z) + O(\log r) \\
    &\leq K_r^{A,B,x}(\|x-y\|_*) + K_r^{A,B,x}(z) + O(\log r).
\end{align*}
Therefore,
\begin{equation*}
    \liminf_{r \to \infty} \frac{K_r^{A,B,x}(y)}{r} \leq \liminf_{r \to \infty} \frac{K_r^{A,B,x}(\|x-y\|_*)}{r} + \limsup_{r \to \infty} \frac{K_r^{A,B,x}(z)}{r} + \limsup_{r \to \infty} \frac{O(\log r)}{r}
\end{equation*}
and so
\begin{align*}
    \dim^{A,B,x}(y) &\leq \dim^{A,B,x}(\|x-y\|_*) + \Dim^{A,B,x}(z) \\
    &\leq\dim^{A,B,x}(\|x-y\|_*) + d - 1.
\end{align*}
In particular, for every $\eps > 0$, the point-to-set principle gives a point $y \in E$ such that
\begin{equation*}
    \hdim E \leq \dim^{A,B,x}(y) + \eps.
\end{equation*}
Thus, for every $\eps$, there exists a point $t=\|x-y\|_* \in \Delta_x^*(E)$ such that
\begin{equation*}
    \hdim E - (d-1) - \eps \leq \dim^{A,B,x}(t).
\end{equation*}
Taking the minimum over all oracles and letting $\eps \to 0$ then gives the desired inequality. \textqed

\section{Sharp examples} \label{s:examples}

In this section we construct the set $E \subset \R^d$ proving Theorem \ref{thm:example}. Fundamentally, the construction adapts \cite{falconer2004fractal} Example 7.8, which in turn works on the same principles as the ``Venetian blinds" constructions that have been studied widely in geometric measure theory for the extreme features they exhibit, e.g., in the context of orthogonal projections, see \cite{csornyei2000visibility}. As such, the use of complexity in \textit{building} our set $E$ is secondary to its use in actually \textit{proving} that it has the asserted properties.

\noindent \define{Proof of Theorem \ref{thm:example}.} The case $s = d$ is trivial, so let $s \in [d-1,d)$ and $\alpha := s - (d-1)$. Since $\| \cdot \|_P$ is a polyhedral norm, there exists a finite set $\{ v^1, ..., v^N \} \subset \R^d$ of vectors\textemdash one for each face of the $\| \cdot \|_P$-ball\textemdash such that 
\begin{equation} \label{T}
    \|x\|_P = \max \, \{ |x \cdot v^\ell| \}_{\ell=1}^N
\end{equation}
for all $x \in \R^d$. Let $c$ be a large constant to be specified later, and let $\{m_k\}_{k=1}^\infty$ be a sequence of positive integers satisfying $m_1 = 1$, $1 + \tfrac{2c}{1-\alpha} < m_2$, and, for all $k \in \Z_+$,
\begin{equation} \label{F}
    k \+ m_k \leq m_{k+1}.
\end{equation}
We also define a sequence $\{n_k\}_{k=1}^\infty$ by
\begin{equation*}
    n_k := m_k + \alpha (m_{k+1} - m_k).
\end{equation*}
Note that, since $\alpha \in [0,1)$ and $m_2 > 1 + \tfrac{2c}{1-\alpha}$, we have $n_k + c < m_{k+1} - c$ for all $k$.

We construct $E \subset [0,1]^d$ in blocks of $N$ steps. Define
\begin{equation} \label{eq:F_k}
    \begin{array}{rl}
        F_k := & \Big\{ \ x \in [0,1]^d :~ \lfloor 2^j( x\cdot v^\ell) \rfloor = 0  \mbox{ for all } \\
        & \ j \in (n_k+c, m_{k+1}-c], \mbox{ where } \ \ell \equiv k \text{ (mod $N$) } \Big\}
    \end{array}
\end{equation}
and
\begin{equation*}
    E := \bigcap_{k=1}^\infty F_k.
\end{equation*} 
These definitions prescribe that the binary expansions of the projections $x \cdot v^\ell$ of the points in $x \in E$ alternate between long strings of free digits and long strings of zeros, with the lengths of the strings rapidly approaching infinity. When a real number has two different binary expansions, we associate to it the expansion terminating in an infinite string of ones. This makes the sets $F_k$ and thus $E$ to be closed.

\noindent \textbf{Claim:} For every oracle $A \subseteq \N$, there exists $x \in E$ such that
\begin{equation} \label{M}
   \liminf_{r \to \infty} \frac{K_r^A(x)}{r} \geq s.  
\end{equation}
Moreover, for all $\tau \in \Delta^P(E)$,
\begin{equation} \label{N}
    \liminf_{r \to \infty} \frac{K_r^A(\tau)}{r} \leq \alpha.
\end{equation}
By the definition of $\alpha$, this will imply that $\hdim \Delta^P(E) \leq s - (d-1) \leq \hdim E - (d-1)$, and an application of Theorem \ref{thm:bound} will entail that these hold with equality.

To obtain this point $x$ and ensure that it belongs to $E$, we iteratively select the digits in the binary expansions of its coordinates. In particular, assuming we have already chosen the digits in the places before the $m_k$th place, the digits in the places $[m_k,m_{k+1})$ will be chosen in two steps. For each $\ell = 1, ..., N$, fix any index $1 \leq i_\ell \leq d$ such that the $i_\ell$th coordinate of $v^\ell$ is nonzero: $v_{i_\ell}^\ell \neq 0$. Let $\ell \equiv k \pmod{N}$, $1 \leq \ell \leq d$.
\begin{enumerate}[label=\arabic*.]
\item We define the digits in the binary expansion of $x_i$ in the places $[m_k,m_{k+1})$ for $i \neq i_\ell$, and in the places $[m_k,n_k)$ for $i = i_\ell$, such that the corresponding $d$ strings (of which $d-1$ have length $m_{k+1} - m_k$ and one has length $n_k - m_k$) are random with respect to each other, the oracle $A$ and the first $m_k$ digits of $x$. In particular, this implies that
\begin{equation*}
    K_{r,m_k}^A(x) \geq d(r-m_k) - O(\log r), \text{ for all }  r \in [m_k,n_k).
\end{equation*}
After we choose the digits of $x_{i_\ell}$ in the places $[n_k,m_{k+1})$ in Step 2, these random choices will allow us to conclude that
\begin{equation*}
    K_{r,n_k}^A(x) \geq (d-\alpha)(r-n_k) - O(\log r), \text{ for all } r \in [n_k,m_{k+1}).
\end{equation*}
\item We want to define the digits in the binary expansion of $x_{i_\ell}$ in the places $[n_k, m_{k+1})$ such that 
\begin{equation} \label{Zo}
    \left\lfloor 2^j \sum_{i=1}^d 2^{-m_{k+1}} \lfloor 2^{m_{k+1}}x_i \rfloor v_i^\ell \right\rfloor \equiv
    \begin{cases}
      0 & \text{if } j \in (n_k+c, m_{k+1}-c],\\
      1 & \text{if } j = m_{k+1}-c+1,\\
      0 & \text{if } j = m_{k+1}-c+2,
    \end{cases}
\end{equation}
where the equivalence holds modulo $2$. Before showing that this is possible, we prove that Equation \eqref{Zo} implies
\begin{equation} \label{Zo2}
    \big\lfloor 2^j (x \cdot v^\ell)\big\rfloor \equiv 0 \pmod{2} \text{ for every } j \in (n_k+c,m_{k+1}-c],
\end{equation}
provided $c$ is chosen large enough. It will then follow from \eqref{eq:F_k} that the point $x$ satisfying \eqref{Zo} belongs to $F_k$. In order to derive Equation \eqref{Zo2} from Equation \eqref{Zo}, we first observe that 
\begin{equation} \label{zoonu}
  x \cdot v^\ell = \sum_{i=1}^d 2^{-m_{k+1}} \lfloor 2^{m_{k+1}}x_i \rfloor v_i^\ell + \sum_{i=1}^d \big( x_i- 2^{-m_{k+1}} \lfloor 2^{m_{k+1}}x_i \rfloor \big) v_i^\ell.
\end{equation}
Notice that, if \eqref{Zo} holds, then for some $M \in \Z_+$ and $t \leq \tfrac{1}{4}$, we have 
\begin{equation} \label{zooni}
    \sum_{i=1}^d2^{-m_{k+1}} \lfloor 2^{m_{k+1}}x_i \rfloor v_i^\ell=2^{-n_{k}-c}M+2^{-m_{k+1}+c} \big( \tfrac{1}{2} + t \big).
\end{equation}
Choose $c$ large enough such that 
\begin{equation} \label{Zoon}
   \max \left( \sum_{i=1}^d |v_i^\ell|, \, \big| v_{i_{\ell}}^{\ell} \big|^{-1} \right) \leq 2^{c-3} \quad \text{for all } \ell,
\end{equation}
from which it follows that
\begin{equation} \label{moon}
\begin{aligned}
    \left| \sum_{i=1}^d(x_i- 2^{-m_{k+1}} \lfloor 2^{m_{k+1}}x_i \rfloor) v_i^\ell \right| & \leq \sum_{i=1}^d \big| x_i- 2^{-m_{k+1}} \lfloor 2^{m_{k+1}}x_i \rfloor \big| |v_i^\ell| \\
    & \leq 2^{-m_{k+1}} \sum_{i=1}^d |v_i^\ell| \\ 
    & \leq 2^{-m_{k+1}+c-3}.
\end{aligned}
\end{equation}
Combining \eqref{zoonu}, \eqref{zooni}, and \eqref{moon}, we see that
\begin{equation*}
    2^j(x \cdot v^\ell) = 2^{j-n_k-c}M + 2^{j-(m_{k+1}-c)} \big( \tfrac{1}{2} + t + t' \big),
\end{equation*}
where $|t'| \leq \tfrac{1}{8}$. Hence, for $j \in (n_{k}+c, m_{k+1}-c]$, the second term above lies in the interval $(0,1)$:
\begin{equation*}
    0 < 2^{j-(m_{k+1}-c)} \big( \tfrac{1}{2} + t + t' \big) < 1.
\end{equation*}
Therefore, we have shown that Equation \eqref{Zo} implies 
\begin{equation} \label{Zoo}
    \left\lfloor 2^j \sum_{i=1}^d 2^{-m_{k+1}} \lfloor 2^{m_{k+1}} x_i \rfloor v_i^\ell \right\rfloor = \big\lfloor 2^j (x \cdot v^\ell) \big\rfloor
\end{equation}
for $j \in (n_k+c, m_{k+1}-c]$, which immediately entails Equation \eqref{Zo2}.

It remains to see that we can choose the digits of $x_{i_\ell}$ in the places $[n_k,m_{k+1})$ such that \eqref{Zo} holds. The condition \eqref{Zo} determines the digits of $x_{i_\ell} \cdot v_{i_\ell}^\ell$ in places $[n_k,m_{k+1}-c+2)$, and we see from \eqref{Zoon} that
\begin{equation*}
    2^{-c+3} \leq v_{i_\ell}^\ell \leq 2^{c-3}.
\end{equation*}
Therefore, we can indeed make such a choice for the digits of $x_{i_\ell}$.
\end{enumerate}
Repeating this process for all $k \in \Z_+$ produces the desired point $x \in E$.

Now we prove that $\dim^A(x) \geq s$. Let $r \in (m_k,m_{k+1}]$ and $\ell \equiv k \pmod N$. We estimate $K_r^A(x)$ in two cases.

\noindent \textbf{Case 1:} $r \in (m_k,n_k]$ \\
By symmetry of information, we have
\begin{equation} \label{A}
\begin{aligned}
    K_r^A(x) & = K_{r,m_k}^A(x) + K_{m_k}^A(x) + O(\log r) \\
    &\geq (r-m_k)d + K_{m_k}^A(x) - O(\log r).
\end{aligned}
 \end{equation}
The inequality is true because the digits of $x$ in the $m_k$th through $r$th places are chosen randomly with respect to $A$ and the digits in the places $1$ through $m_k$.

\noindent \textbf{Case 2:} $r \in (n_k,m_{k+1}]$ \\
By the symmetry of information, we have
\begin{equation*}
    K_r^A(x) = K_{r,n_k}^A(x) + K_{n_k}^A(x) - O(\log r).
\end{equation*}
It follows from \eqref{A} of Case 1 that 
\begin{align*}
    K_{n_k}^A(x) & \geq (n_k-m_k)d + K_{m_k}^A(x) - O(\log r) \\
    & = \alpha(m_{k+1}-m_k)d + K_{m_k}^A(x) - O(\log r).
\end{align*}
Thus we have 
\begin{equation} \label{E}
\begin{aligned}
    K_r^A(x) & \geq K_{r,n_k}^A(x) + \alpha(m_{k+1}-m_k)d+K_{m_k}^A(x) - O(\log r) \\
    & \geq K_{r,n_k}^A(x_1,...,x_{i_{\ell}-1},x_{i_{\ell}+1},...,x_d) + \alpha(m_{k+1}-m_k)d + K_{m_k}^A(x) - O(\log r) \\
    & \geq (r-n_k)(d-1) + \alpha(m_{k+1}-m_k)d + K_{m_k}^A(x) - O(\log r) \\
    & = (r-m_k)(d-1) + \alpha(m_{k+1}-m_k) + K_{m_k}^A(x) - O(\log r).
\end{aligned}
\end{equation}
The third inequality above also follows from randomness. We put $r = m_{k+1}$ in \eqref{E} to get
\begin{align*}
    K_{m_{k+1}}^A(x) &= (m_{k+1}-m_k)s + K_{m_k}^A(x) - O(\log r) \\
    &\geq m_{k+1} \left( 1 - \frac{1}{k} \right) s + K_{m_k}^A(x) - O(\log r),
\end{align*}
where the inequality follows from \eqref{F}. Thus for any arbitrarily small $t > 0$, we have 
\begin{equation*}
   \frac{K_{m_{k+1}}^A(x)}{m_{k+1}} \geq s-t \quad  \text{for large enough } k.
\end{equation*}
We will now justify that the expression $K_r^A(x)/r$ is minimized when $r = m_k$ for some $k$. By \eqref{A} of Case 1,
\begin{equation} \label{x}
\begin{aligned}
    \frac{K_r^A(x)}{r} & \geq \left( 1-\frac{m_k}{r} \right) d + \frac{K_{m_k}^A(x)}{r} - \frac{O(\log r)}{r} \\
    & = \left( 1-\frac{m_k}{r} \right) d + \frac{K_{m_k}^A(x)}{m_{k}} \left( \frac{m_k}{r} \right) - \frac{O(\log r)}{r} \\
    & \geq \left( 1-\frac{m_k}{r} \right) d + (s-t) \frac{m_k}{r} - \frac{O(\log r)}{r} \\
    & = d - (\alpha-1+t) \frac{m_k}{r} - \frac{O(\log r)}{r}
\end{aligned}
\end{equation}
and similarly, by \eqref{E}, we have in Case 2
\begin{equation} \label{y}
\begin{aligned}
    \frac{K_r^A(x)}{r} &\geq \left( 1-\frac{m_k}{r} \right)(d-1) + \frac{\alpha(m_{k+1}-m_k)}{r} +\frac{K_{m_k}^A(x)}{r} - \frac{O(\log r)}{r} \\
    &\geq d-1 + \alpha\frac{m_{k+1}}{r} + \frac{K_{m_k}^A(x) - sm_k}{r}   - \frac{O(\log r)}{r} \\
    & \geq d - 1 + \alpha\frac{m_{k+1}}{r} - \frac{tm_k}{r} - \frac{O(\log r)}{r}.
\end{aligned}
\end{equation}
Recalling $\tfrac{m_k}{r} \leq 1$ and $\tfrac{m_{k+1}}{r} \geq 1$, we conclude \eqref{M} from \eqref{x} and \eqref{y}, and therefore we have proved that $\hdim E \geq s$.

Now, for any fixed $z,y \in E$, let $\tau := \| z-y \|_P \in \Delta^P(E)$. Thus, $\tau = |(z-y) \cdot v^\ell|$ for some $\ell$. Hence, for any $k \equiv \ell \pmod N$, we have
\begin{equation*}
    \lfloor 2^j((z-y) \cdot v^\ell) \rfloor \equiv 0 \pmod 2 \quad \text{for all} \quad j \in [n_k+c+1, m_{k+1}-c-1]
\end{equation*}
or
\begin{equation*}
    \lfloor 2^j((z-y) \cdot v^\ell) \rfloor \equiv 1 \pmod 2 \quad \text{for all} \quad j \in [n_k+c+1, m_{k+1}-c-1].
\end{equation*}
Therefore,
\begin{equation*}
    K_{m_{k+1}-c,n_k+c}^A(\tau) \leq O(\log m_{k+1}).
\end{equation*}
This is because specifying a string of 000's or 111's having length less than $m_{k+1}$ requires a program of length at most $O(\log m_{k+1})$. We consequently have 
\begin{align*}
    K_{m_{k+1}-c}^A(\tau) & \leq K_{m_{k+1}-c,n_k+c}^A(\tau) + K_{n_k}^A(\tau) + O(\log m_{k+1}) \\
    & \leq n_k + O(\log m_{k+1}) \\
    & \leq m_{k+1}\alpha \left( 1-\frac{1}{k} \right) + O(\log m_{k+1}),
\end{align*}
where the last inequality follows from \eqref{F}. Taking $k \to \infty$ gives \eqref{N}, and, since this holds for all $\tau \in \Delta^P(E)$, we conclude that $\hdim \Delta^P(E) \leq \alpha$. The desired upper bound on $\hdim E$ then follows immediately from Theorem \ref{thm:bound}:
\begin{equation*}
    s \leq \hdim E \leq \hdim \Delta^P(E) + (d-1) \leq \alpha + (d-1) = s.
\end{equation*}
As such, all the inequalities above are actually equalities, so Equation \eqref{eq:sharp} holds. \textqed

\section*{Acknowledgement}
The authors thank Kaiyi Huang for her careful reading and insights concerning Theorem \ref{thm:example}. We would also like to thank the editor and an anonymous referee for many helpful corrections and suggestions.

\bibliographystyle{plain}
\bibliography{references}

\newpage

\end{document}